\pdfoutput=1
\documentclass[11pt, oneside]{article}   	
\usepackage{geometry}                		
\usepackage{graphicx}				
\usepackage{caption}								
\usepackage{subcaption}
\usepackage{tabularx}
\usepackage{booktabs}

\usepackage{amssymb}
\usepackage{url}
\usepackage{hyperref}
\usepackage{amsmath}

\title{Periodic modification of the Boerdijk-Coxeter helix (tetrahelix)}
\author{Garrett Sadler\footnote{Corresponding author: \href{mailto:garrett@quantumgravityresearch.org}{garrett@quantumgravityresearch.org}}, Fang Fang, Julio Kovacs, Klee Irwin\\\emph{Quantum Gravity Research}, Topanga, CA, USA}

\begin{document}
\maketitle

\begin{abstract}
The Boerdijk-Coxeter helix \cite{coxeter1974}\cite{boerdijk1952} is a helical structure of tetrahedra which possesses no non-trivial translational or rotational symmetries. In this document, we develop a procedure by which this structure is modified to obtain both translational and rotational (upon projection) symmetries along/about its central axis. We report the finding of several, distinct periodic structures, and focus on two particular forms related to the pentagonal and icosahedral aggregates of tetrahedra as well as Buckminster Fuller's ``jitterbug transformation''.
\end{abstract}

\section{Introduction}\label{intro}

\begin{figure}
	\centering
	\begin{subfigure}[t]{0.3\textwidth}
		\centering
		\includegraphics[width=\textwidth]{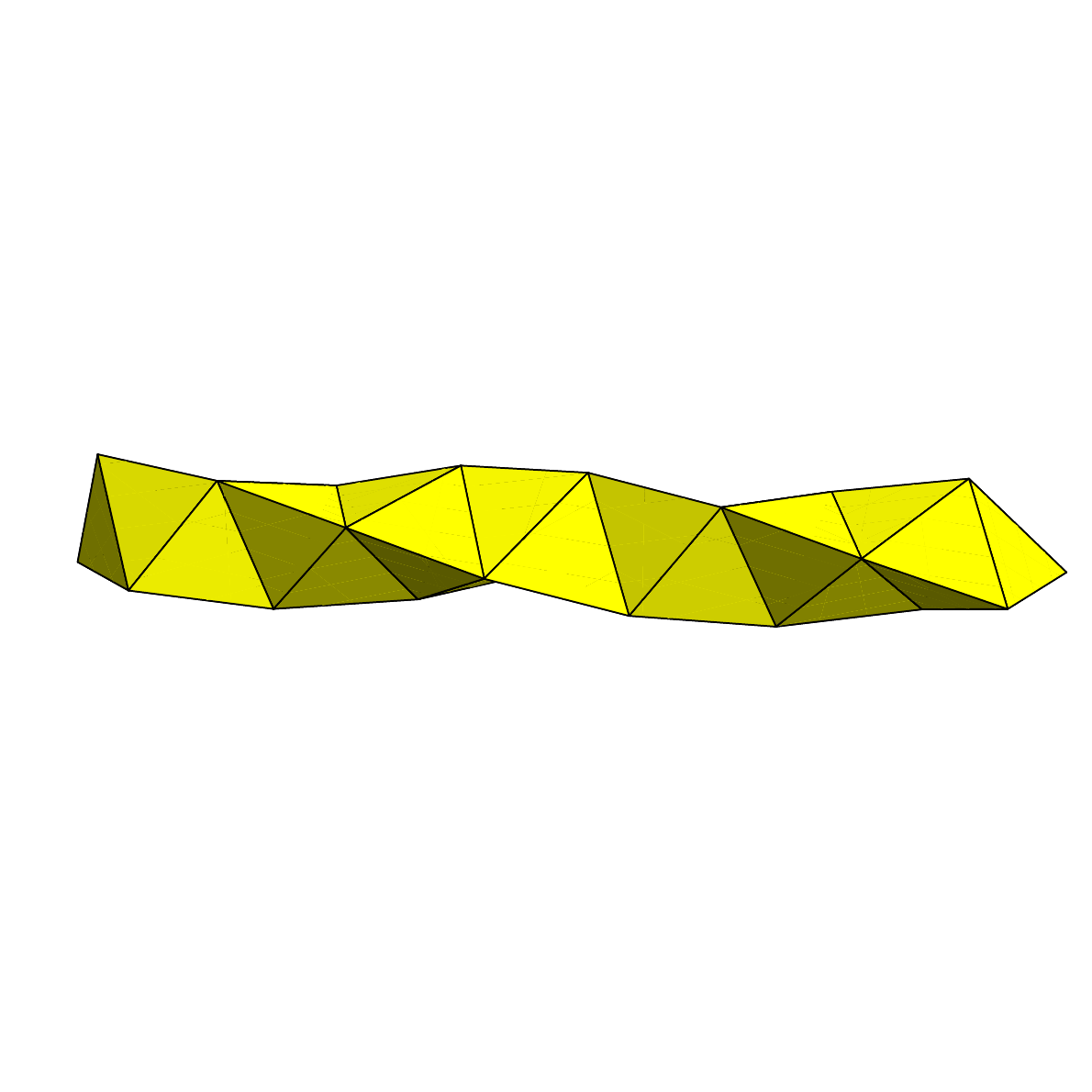}
		\caption{A right-handed BC helix.}
		\label{F:philicesA}
	\end{subfigure}
	\quad
	\begin{subfigure}[t]{0.3\textwidth}
		\centering
		\includegraphics[width=\textwidth]{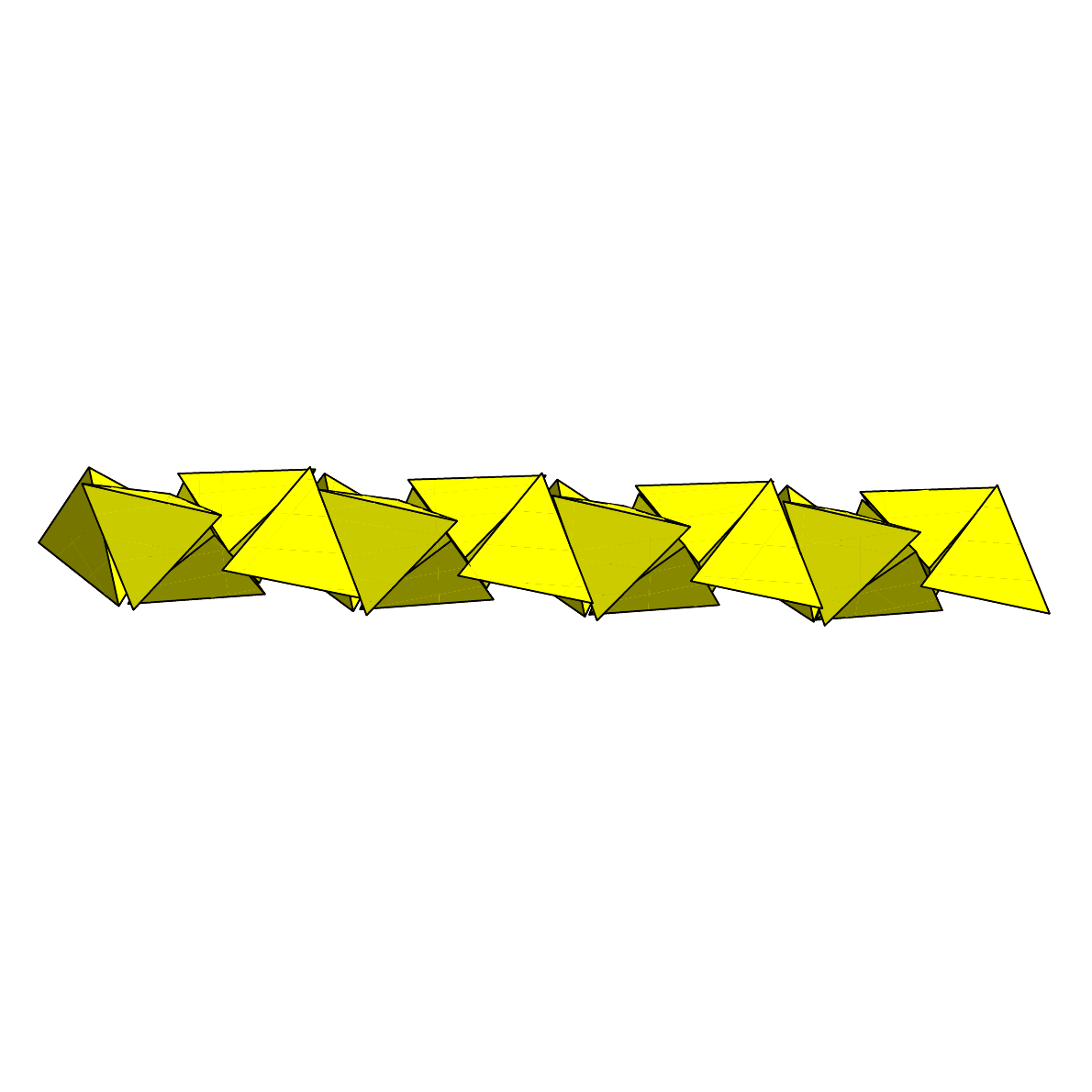}
		\caption{A ``5-BC helix'' may be obtained by appending and rotating tetrahedra through the angle given by \eqref{E:beta} using the same chirality of the underlying helix.}
		\label{F:philicesB}
	\end{subfigure}
	\quad
	\begin{subfigure}[t]{0.3\textwidth}
		\centering
		\includegraphics[width=\textwidth]{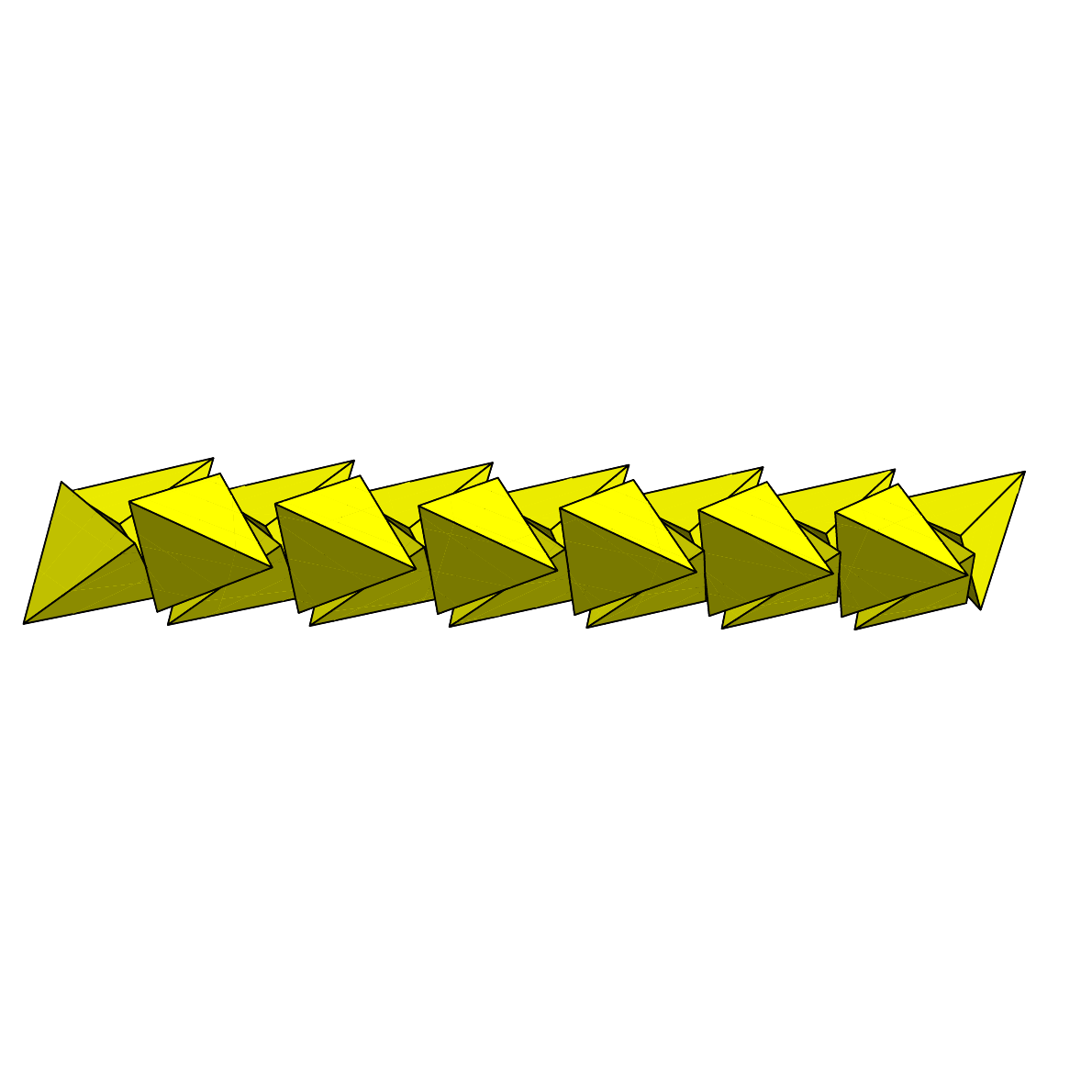}
		\caption{A ``3-BC helix'' may be obtained by appending and rotating tetrahedra through the angle given by \eqref{E:beta} using the opposite chirality of the underlying helix.}
		\label{F:philicesC}
	\end{subfigure}
	\caption{Canonical and modified Boerdijk-Coxeter helices.}\label{F1:philices}
\end{figure}

The Boerdijk-Coxeter helix (BC helix, tetrahelix) \cite{coxeter1974}\cite{boerdijk1952} is an assemblage of regular tetrahedra in a linear, helical fashion (Figure~\ref{F:philicesA}). This assemblage may be obtained by appending faces of tetrahedra together so as to maintain a central axis or, alternatively, R.W. Gray \cite{gray} has produced a description of the BC helix by partitioning into 4-tuples the points of $\mathbb{R}^3$ given by the sequence $\left(s_n\right)_{n \in \mathbb{Z}}$
\newline
\begin{equation}
	s_{n} = \left( r \cos{( n \theta ) } , \pm r \sin{(n \theta)} , n h \right),
\end{equation}\label{parametrictetrahelix}
\newline
where $r = 3 a \sqrt{3} / 10 $, $\theta = \arccos (- 2/3 )$, $h = a / \sqrt{10}$, and $a$ designates the tetrahedral edge length. (The sequence of faces used while appending, or the sign of the second term in \eqref{parametrictetrahelix}, determine the chirality of the helix.) Due to the irrational value of $\theta$, it may be observed that the BC helix has an aperiodic nature, in that the structure has no non-trivial translational or rotational symmetries. Here, we describe a modified form of the BC helix that has both translational and rotational symmetries along/about its central axis. Figures~\ref{F:philicesB} and \ref{F:philicesC} show two such modified structures.

\section{Method of assembly: modified BC helices}\label{S:genphilix}
\begin{figure}
	\centering
	\begin{subfigure}[t]{0.3\textwidth}
		\centering
		\includegraphics[width=\textwidth]{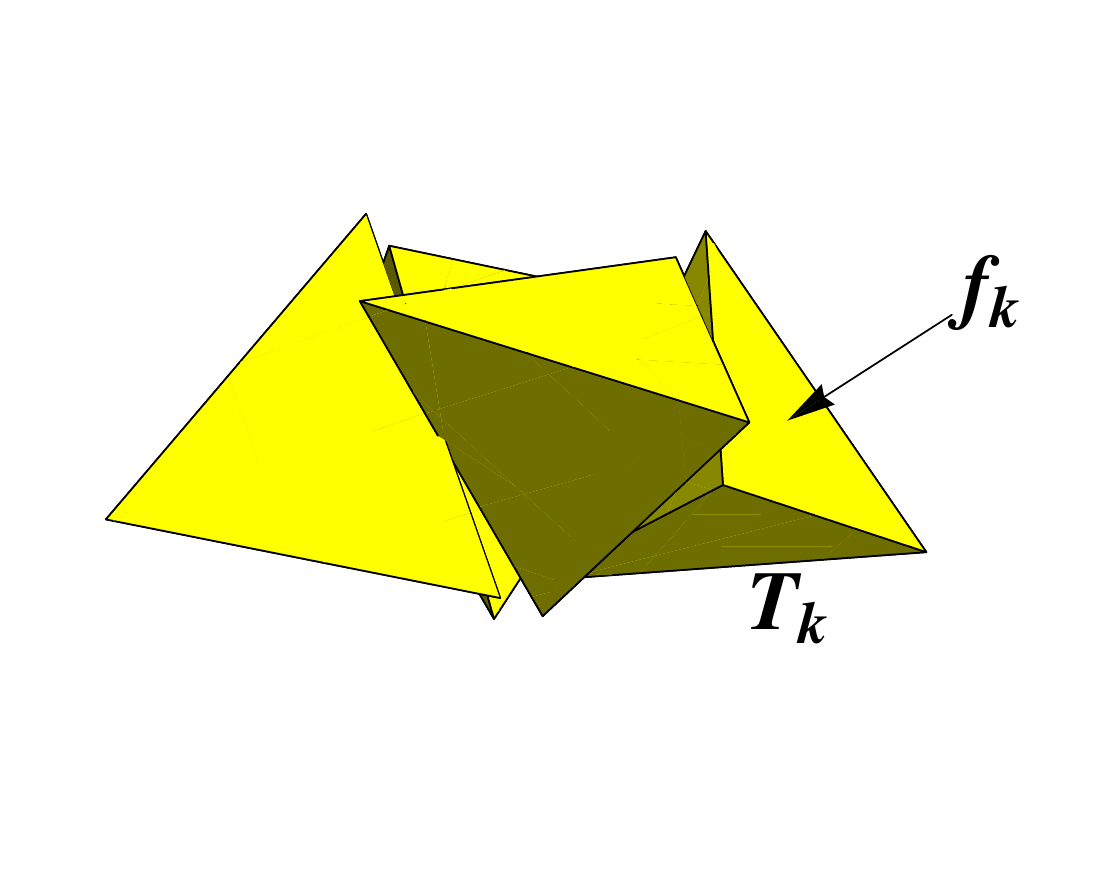}
		\caption{A segment of an $m$-BC helix with face $f$ identified on tetrahedron $T_k$.}
		\label{F:makeaphilixA}
	\end{subfigure}
	\quad
	\begin{subfigure}[t]{0.3\textwidth}
		\centering
		\includegraphics[width=\textwidth]{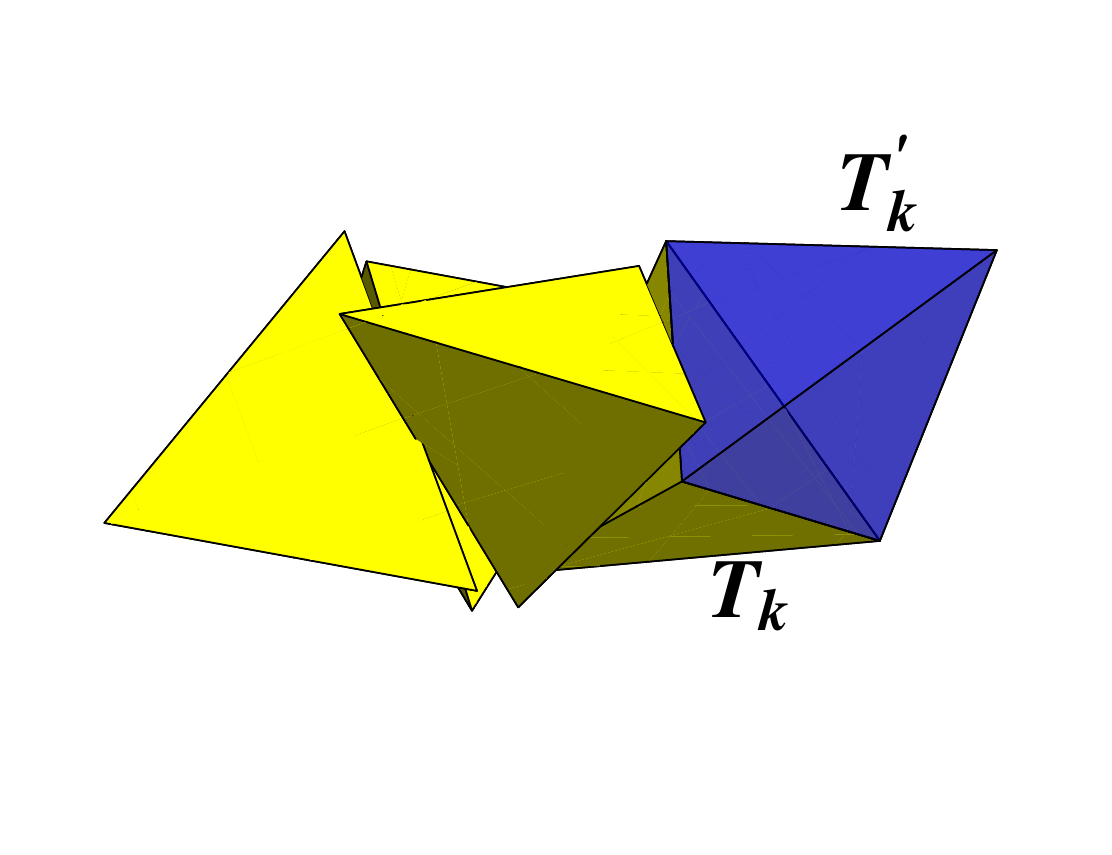}
		\caption{An interim tetrahedron, $T_k^\prime$ (shown in blue), is appended (face-to-face) to $f$ on $T_k$.}
		\label{F:makeaphilixB}
	\end{subfigure}
	\quad
	\begin{subfigure}[t]{0.3\textwidth}
		\centering
		\includegraphics[width=\textwidth]{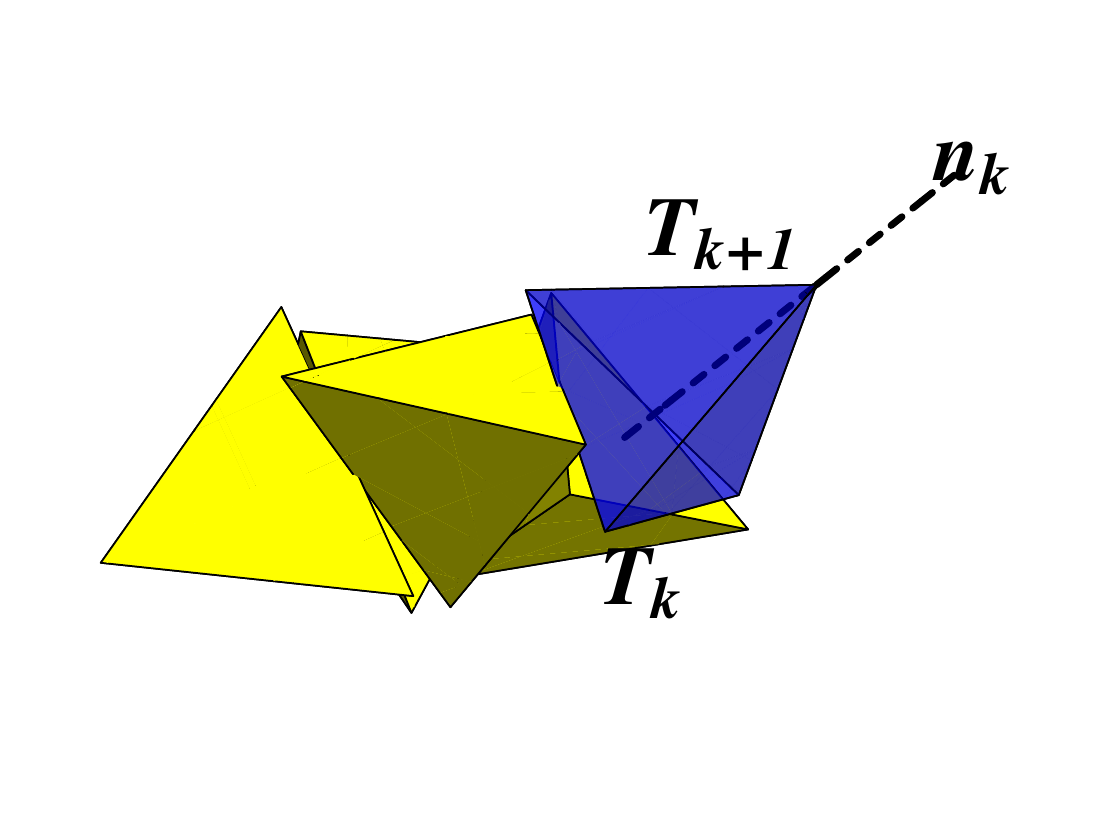}
		\caption{Finally, $T_{k+1}$ is obtained by rotating $T_k^\prime$ through the angle $\beta$ about the axis $n_k$.}
		\label{F:makeaphilixC}
	\end{subfigure}
	\caption{Assembly of modified BC helix.}\label{F:makeaphilix}
\end{figure}

The assembly our modified BC helices is distinguished from that of the canonical BC helix in that an additional operation is required between appending tetrahedra to the helix. This operation is depicted in Figure~\ref{F:makeaphilix}. Starting with a tetrahedron $T_k = \left( v_{k 0}, v_{k 1}, v_{k 2}, v_{k 3} \right)$, a face $f_k$ is selected onto which an interim tetrahedron, $T_k^\prime$, is appended. The $\left( k+1 \right)^{\text{th}}$ tetrahedron is obtained by rotating $T_k^\prime$ through an angle $\beta$ about an axis $n_k$ normal to $f_k$, passing through the centroid of $T_k^\prime$. 

The resulting structure depends, principally, on two choices in this process. Firstly, as with the BC helix, the sequence of faces, $F = \left(f_0, f_1,\ldots,f_k \right)$, selected in the construction of the helix will determine its \emph{underlying chirality}---i.e., the chirality of the helix formed by the tetrahedral centroids. Faces may be selected so that some sequences produce right-handed helices, while others produce left-handed helices (and, certainly, some sequences do not produce helices at all). Secondly, there is the choice of the magnitude and direction of the rotation. In the present writing, we will use the convention that a facial normal vector $n_k$ is pointed away from the face $f_k$, i.e., $n_k$ points away from the interior of $T_k$. Consequently, positive values of $\beta$ will correspond to right-handed rotations about $n_k$, while negative values will produce left-handed rotations. (And, certainly, a canonical BC helix is obtained for $\beta = 0$.)

A convenient method of assembly for a modified BC helix is by usage of two transformations, $A_T^f\colon \mathbb{R}^3 \longrightarrow \mathbb{R}^3$ and $B_T^f\colon \mathbb{R}^3 \longrightarrow \mathbb{R}^3$, where $A_T^f$ is a reflection across face $f$ on tetrahedron $T$, and $B_T^f$ is a rotation about an axis normal to this face, passing through its center.

A modified BC helix is formed by applying $A_{T_k}^{f_k}$ to the vertices $v_{k j},\ j = 0,\ldots,3$, of tetrahedron $T_k$ to produce $T_k^\prime$. Finally, $T_{k+1}$ is obtained by applying $B_{T_k^\prime}^{f_k}$ to the vertices of $T_k^\prime$, that is:
\newline
\begin{align}
	v_{k j}^\prime &= A_{T_k}^{f_k} \left( v_{k j} \right)\\
	v_{\left(k+1\right) j} &= B_{T_k^\prime}^{f_k} \left( v_{k j}^\prime \right).
\end{align}
\newline
By applying these transformations in an alternating fashion, first to each $T_k$ and then to $T_k^\prime$, a modified BC helix is assembled.

When referring to a modified BC helix, we use the term \emph{period} to refer to the number of appended tetrahedra necessary to return to an initial angular position on the helix, and will say that the structure is \emph{periodic} when such an integer exists. For almost all values of $\beta$, the associated modified BC helix is aperiodic, however, the resulting structure is periodic for certain values of $\beta$. Here, we use the term \emph{m-BC~helix} to designate a modified BC helix with a period of $m$ tetrahedra.

Using analytical and numerical methods, our group has found modified BC helices with periods of $2,\ldots,20$ tetrahedra, with the exception of a 6-BC helix. In an upcoming paper, we demonstrate the impossibility of a modified BC helix with a period of 6 tetrahedra. Table~\ref{T:periods} provides a summary of these BC helix periodicities and their corresponding values of $\beta$, assuming an assembly face sequence $F$ that produces a right-handed helix. (For a sequence that produces a left-handed helix, one only needs to change the sign of the angular values in Table~\ref{T:periods}.)

\begin{table}[t]
	\centering
	\begin{tabular}{  c c c c c c }
		\toprule
		\multicolumn{2}{c}{Exact values} & \multicolumn{4}{c}{Numerically obtained values}\\
		\cmidrule(r){1-2}\cmidrule(l){3-6}
		Period & $\beta$ & Period & $\beta$ & Period & $\beta$ \\
		\midrule
		2	& $\frac{\pi}{3}$						& 7	& $-0.69115$	& 14 & $-1.51006$ \\
			&								&	& $0.494277$	&	& $-0.0733038$\\
		3	& $-\arccos \left( \frac{3 \phi - 1}{4} \right)$	& 8	& $0.0712094$	& 15	& $-0.873363$ \\
			&								&	&			&	& $0.789587$ \\
		4	& $-\frac{\pi}{3}$					& 9	& $-1.38858$	& 16 & $-0.450295$ \\
			&								&	& $0.617847$	&	& $0.565487$ \\
		5	& $\arccos \left( \frac{3 \phi - 1}{4} \right)$	& 10	& $-0.559203$	& 17	& $-0.613658$ \\
			&								&	&			&	& $0.821003$ \\
			&								& 11 & $-0.81472$	& 18 & $-0.764454$ \\
			&								&	& $0.697434$	&	& $0.182212$ \\
			&								& 12	& $0.402124$	& 19	& $-0.908967$ \\
			&								&	&			&	& $0.844041$ \\
			&								& 13	& $-0.492183$	& 20	& $-0.131947$ \\
			&								&	& $0.751888$	&	& $0.661829$ \\
		\bottomrule
	\end{tabular}
	\caption{\emph{m}-BC helix periods and their corresponding values of $\beta$, assuming tetrahedra are appended in a sequence $F$ which produces a right-handed helix. For values of $\beta$ corresponding to a left-handed underlying helix, the sign of $\beta$ must be changed. In this table, $\phi$ designates the golden ratio.}
	\label{T:periods}
\end{table}

In \cite{goldenQC}, we present novel modifications of icosahedral and pentagonal bipyramid aggregates of tetrahedra involving a rotation through an angular value of
\newline
\begin{equation}\label{E:beta}
	\beta = \pm \arccos \left( \frac{3 \phi - 1}{4} \right),
\end{equation}
\newline
where $\phi = \left( 1 + \sqrt{5} \right)/2$ denotes the golden ratio. As indicated in Table~\ref{T:periods}, this value of $\beta$ corresponds to 3-- and 5-BC helices. For this reason, as well as the appearance of this angle in Fuller's ``jitterbug transformation''  \cite{fuller1975}, we will focus on the 3-- and 5-BC helix in the present writing. In Section~\ref{S:fperiod}, we will provide an explicit construction of a 5-BC helix, along with some additional properties of this structure. In Section~\ref{SS:3period}, the same is done for the 3-BC helix.

\section{Modified BC helices: explicit examples}\label{S:simplephilix}

In this section we will describe the assembly of the 3-- and 5-BC helices. The approach used here generates a primitive set of tetrahedra following the method of assembly described in Section~\ref{S:genphilix} while using the value of $\beta$ in \eqref{E:beta}. Modified BC helices of arbitrary length may then be generated by translating copies of this primitive set along the helix's central axis (explicitly provided below). Due to the presence of the golden ratio in \eqref{E:beta}, we refer to such a structure by the name ``\emph{phi}lix''.

In order to keep the expressions simple, we choose the starting tetrahedron in a convenient way. The expressions for any desired philix axis can be obtained by multiplying the values given here by the corresponding rotation matrix. At the conclusion of each of the sections below, the appropriate transformation is offered to align the philix axis with the \emph{z}-axis of $\mathbb{R}^3$.

Interestingly, the sign of $\beta$ will determine whether a 3-- or a 5-period philix is generated according to the following rule:

\begin{quote}
	(\emph{i}) \emph{When the chiralities of the rotation by} $\beta$ \emph{and that of the underlying helix produced by the face sequence} $F=\left(f_0,f_1,\ldots,f_k\right)$ \emph{are} \textbf{alike}\emph{, one obtains a 5-period philix.}
	\newline
	\newline
	(\emph{ii}) \emph{When the chiralities of the rotation by} $\beta$ \emph{and that of the underlying helix produced by the face sequence} $F=\left(f_0,f_1,\ldots,f_k\right)$ \emph{are} \textbf{unlike}\emph{, one obtains a 3-period philix.}
\end{quote}

In the constructions of Sections~\ref{S:fperiod} and \ref{SS:3period}, face sequences are used such that a right-handed underlying helix is produced. Accordingly, a positive value of $\beta$ generates a 5-BC helix, while a negative value generates a 3-BC helix. For compactness, the values of the primitive tetrahedral vertices, central axis vector, and central helix radius and pitch are given in these sections. All values and expressions necessary to compute the transformations $A_{T_k}^{f_k}$ and $B_{T_k}^{f_k}$ are given in Appendix~A.

\subsection{The 5-BC helix}\label{S:fperiod}
Using $T_{k} = \left( v_{k 0}, v_{k 1}, v_{k 2}, v_{k 3} \right)$, $ v_{k j} \in \mathbb{R}^3$, to designate a tetrahedron of an 5-BC helix, a primitive set for a 5-period philix may be formed from the unit-edge length tetrahedra $\{ T_{0},\dots,T_{4} \}$ given by

\allowdisplaybreaks{
\begin{align}
	T_{0}\colon \qquad 	v_{0 0} &= \left( 0, 0, \sqrt{\frac{2}{3}} - \frac{1}{2 \sqrt{6}} \right) \label{T0}\\
					v_{0 1} &= \left( - \frac{1}{2 \sqrt{3}}, - \frac{1}{2}, - \frac{1}{2 \sqrt{6}} \right)\notag\\
					v_{0 2} &= \left( - \frac{1}{2 \sqrt{3}}, \frac{1}{2}, - \frac{1}{2 \sqrt{6}} \right)\notag\\
					v_{0 3} &= \left( \frac{1}{\sqrt{3}}, 0, - \frac{1}{2 \sqrt{6}} \right)\notag\\				
	T_{1}\colon \qquad	v_{1 0} &= \left( 0, 0, - \frac{5}{2 \sqrt{6}} \right) \label{T1}\\
					v_{1 1} &= \left( - \frac{1+3 \sqrt{5}+3 \sqrt{6-2 \sqrt{5}}}{16 \sqrt{3}}, - \frac{1+3 \sqrt{5}- \sqrt{6-2 \sqrt{5}}}{16}, - \frac{1}{2 \sqrt{6}} \right)\notag\\
					v_{1 2} &= \left( - \frac{1+3 \sqrt{5}-3 \sqrt{6-2 \sqrt{5}}}{16 \sqrt{3}}, \frac{1+3 \sqrt{5}+ \sqrt{6-2 \sqrt{5}}}{16}, - \frac{1}{2 \sqrt{6}} \right)\notag\\
					v_{1 3} &= \left( \frac{1+3 \sqrt{5}}{8 \sqrt{3}}, - \frac{1}{4} \sqrt{ \frac{1}{2} \left(3-\sqrt{5} \right)}, - \frac{1}{2 \sqrt{6}} \right)\notag\\
	T_{2}\colon \qquad	v_{2 0} &= \left( - \frac{1}{12 \sqrt{3}}, \frac{-4+ \sqrt{5}}{12}, - \frac{8 + 3\sqrt{5}}{6 \sqrt{6}} \right)\label{T2}\\
					v_{2 1} &= \left( - \frac{11 + 3 \sqrt{5}}{24 \sqrt{3}}, - \frac{5 + \sqrt{5}}{24}, \frac{-8 + 3 \sqrt{5}}{6 \sqrt{6}} \right)\notag\\
					v_{2 2} &= \left( \frac{5 - 3 \sqrt{5}}{12 \sqrt{3}}, \frac{5 + \sqrt{5}}{12}, - \frac{5}{6 \sqrt{6}} \right)\notag\\
					v_{2 3} &= \left( - \frac{5}{72} \left(\sqrt{3} + 3 \sqrt{15} \right), \frac{5}{24} \left(-1 + \sqrt{5} \right), - \frac{11}{6 \sqrt{6}} \right)\notag\\
	T_{3}\colon \qquad	v_{3 0} &= \left( \frac{5 - 4 \sqrt{5}}{12 \sqrt{3}}, - \frac{\sqrt{5}}{12}, - \frac{11 + 2 \sqrt{5}}{6 \sqrt{6}} \right)\label{T3}\\
					v_{3 1} &= \left( \frac{13 - 11 \sqrt{5}}{24 \sqrt{3}}, \frac{ 3 + 7 \sqrt{5}}{24}, - \frac{8 + 5 \sqrt{5}}{6 \sqrt{6}} \right)\notag\\
					v_{3 2} &= \left( \frac{13 - 5 \sqrt{5}}{24 \sqrt{3}}, \frac{ -3 + 7 \sqrt{5}}{24}, \frac{- 8 + 5 \sqrt{5}}{6 \sqrt{6}} \right)\notag\\
					v_{3 3} &= \left( - \frac{5 + 2 \sqrt{5}}{6 \sqrt{3}}, \frac{\sqrt{5}}{6}, - \frac{5 + 2 \sqrt{5}}{6 \sqrt{6}} \right)\notag\\
	T_{4}\colon \qquad	v_{4 0} &= \left( \frac{5 \left(1 - \sqrt{5} \right)}{12 \sqrt{3}}, \frac{-5 + \sqrt{5}}{12}, - \frac{5 + 4 \sqrt{5}}{6 \sqrt{6}} \right)\label{T4}\\
					v_{4 1} &= \left( - \frac{5 + \sqrt{5}}{24 \sqrt{3}}, \frac{5 \left(1 + \sqrt{5} \right)}{24}, - \frac{11 + 4 \sqrt{5}}{6 \sqrt{6}} \right)\notag\\
					v_{4 2} &= \left( - \frac{11 + 13 \sqrt{5}}{24 \sqrt{3}}, \frac{5 - \sqrt{5}}{24}, - \frac{8 + 7 \sqrt{5}}{6 \sqrt{6}} \right)\notag\\
					v_{4 3} &= \left( - \frac{1 + 8 \sqrt{5}}{12 \sqrt{3}}, \frac{4 + \sqrt{5}}{12}, - \frac{8 + \sqrt{5}} {6 \sqrt{6}} \right).\notag	
\end{align}
}
\newline
A 5-period philix may be generated by translating the vertices of these tetrahedra by integer values of a vector $w_5 \in \mathbb{R}^3$ given by
\newline
\begin{equation}\label{axisz}
	w_5 = \left( - \frac{ 5 \left( \sqrt{3} + \sqrt{15} \right)}{36}, \frac{5 + \sqrt{5}}{12}, - \frac{5 + 2 \sqrt{5}}{3 \sqrt{6}} \right),
\end{equation}
\newline
such that
\newline
\begin{equation}\label{periodicrelation}
	v_{\left( j + 5 k \right) i} = v_{j i} + k w_5, \text{\qquad for } k \in \mathbb{Z}.
\end{equation}
\begin{figure}
	\centering
	\begin{subfigure}[t]{0.45\textwidth}
		\centering
		\includegraphics[width=\textwidth]{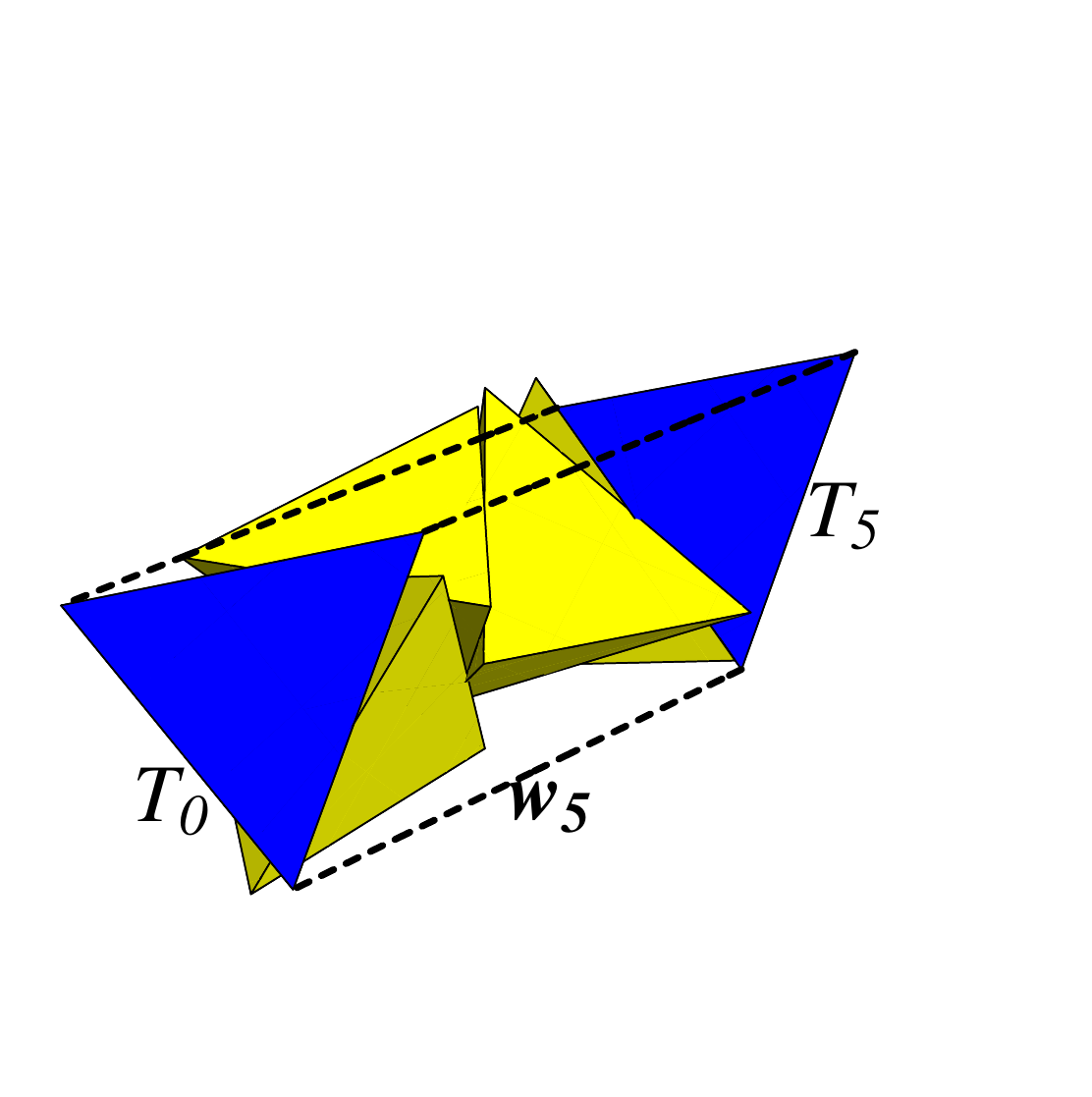}
		\caption{The vertices of $T_5$ are the vertices of $T_0$ translated by $w_5$.}
		\label{F2:5philixA}
	\end{subfigure}
	\quad
	\begin{subfigure}[t]{0.45\textwidth}
		\centering
		\includegraphics[width=\textwidth]{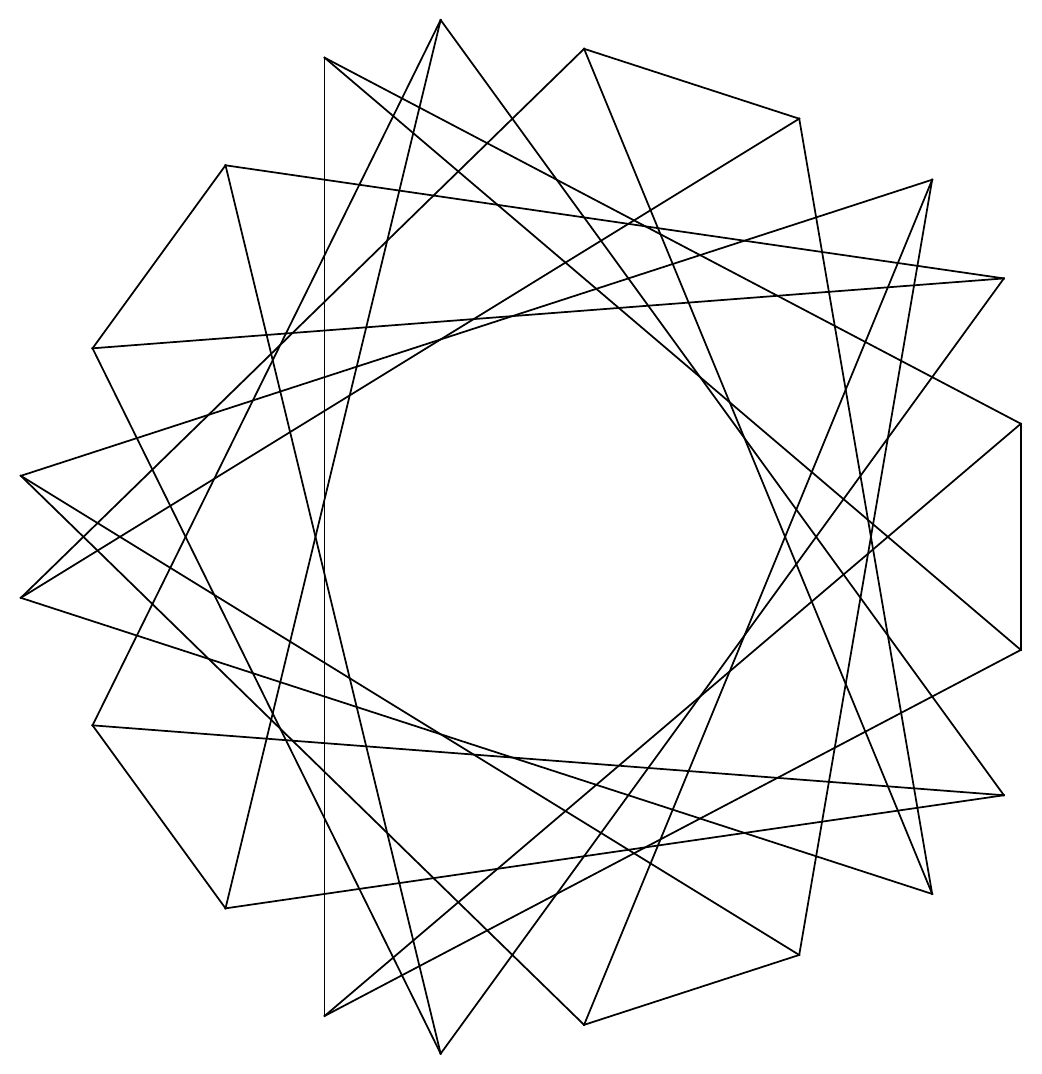}
		\caption{A projection of the 5-BC helix along its central axis.}
		\label{F2:5philixB}
	\end{subfigure}
	\caption{The periodicity of the 5-BC helix.}
	\label{F2:5philix}
\end{figure}

When this is done, one obtains a structure with 5-fold rotational symmetry (in its projection) and a linear ``period'' of 5 tetrahedra along its central axis (see Figure~\ref{F2:5philix}). The centroids of the tetrahedra comprising a 5-period philix form a helix with a linear pitch of
\newline
\begin{equation}\label{5pitch}
	p_5 = \sqrt{\frac{25}{18} + \frac{5 \sqrt{5}}{9}}
\end{equation}
\newline
and a radius of
\newline
\begin{equation}\label{5radius}
	r_5 = \frac{5-\sqrt{5}}{15\sqrt{2}}
\end{equation}
\newline
producing a helix with the parameterization $c\colon \mathbb{R}\longrightarrow\mathbb{R}^3$ given by:
\newline
\begin{equation}\label{5helixeq}
	c \left( t \right) = r_5 \left( u_1 \cos{t} + u_2 \sin{t} \right) + \frac{t}{4\pi} w_5 + q_5,
\end{equation}
\newline
where $u_1 = \left( -\frac{1}{\sqrt{6}}, \frac{1}{\sqrt{2}}, \frac{1}{\sqrt{3}} \right)$ and $u_2 = \left(-\frac{1}{2} \sqrt{\frac{1}{3} \left(5+\sqrt{5}\right)},-\frac{1}{2} \sqrt{1+\frac{1}{\sqrt{5}}},\frac{1}{\sqrt{15+6 \sqrt{5}}}\right)$ are orthonormal vectors spanning the plane perpendicular to the philix axis $w_5$, and 
\newline
\begin{equation}
	q_5=\left(-\frac{\sqrt{5}-5}{30 \sqrt{3}},\frac{1}{30} \left(\sqrt{5}-5\right),\frac{\sqrt{5}-5}{15 \sqrt{6}}\right)
\end{equation}
\newline
is a vector to translate the helix to the location of the philix above (as its axis does not pass through the origin). The tetrahedral centroids lie on this helix at the positions given by $t = k \frac{4 \pi}{5}, k \in \mathbb{Z}$.

The 5-period philix described in this section may by aligned with the \emph{z}-axis by applying the transformation
\newline
\begin{equation}\label{H5transform}
	H_5 \left(v\right) = C_5 \left(v - q_5 \right),
\end{equation}
\newline
to each vertex $v_{k j}$ of the philix, where
\newline
\newcommand\scalemath[2]{\scalebox{#1}{\mbox{\ensuremath{\displaystyle #2}}}}
\begin{equation}
	C_5 = \scalemath{0.8}{
	\left(
\begin{array}{ccc}
 \frac{1}{24} \left(9-\sqrt{75+30 \sqrt{5}}\right) & \frac{1}{4} \sqrt{\frac{1}{2} \left(10+\sqrt{5}+\sqrt{75+30 \sqrt{5}}\right)} & \frac{\left(3+\sqrt{5}\right)
   \left(5+\sqrt{5}\right) \sqrt{6 \left(5+2 \sqrt{5}\right)}}{300+132 \sqrt{5}} \\
 \frac{1}{4} \sqrt{\frac{1}{2} \left(10+\sqrt{5}+\sqrt{75+30 \sqrt{5}}\right)} & \frac{5}{8}-\frac{1}{8} \sqrt{3+\frac{6}{\sqrt{5}}} & -\frac{1}{2}
   \sqrt{1-\frac{1}{\sqrt{5}}} \\
 -\frac{25 \left(123+55 \sqrt{5}\right)}{2 \sqrt{6} \left(5+2 \sqrt{5}\right)^{7/2}} & \frac{1}{2} \sqrt{1-\frac{1}{\sqrt{5}}} & -\frac{5 \left(360+161
   \sqrt{5}\right)}{\sqrt{3} \left(5+2 \sqrt{5}\right)^{7/2}}
\end{array}
\right)
	}
\end{equation}
\newline
is a matrix that rotates $w_5$ to the direction of $\left(0,0,1\right)$.
\begin{figure}
	\centering
	\begin{subfigure}[t]{0.45\textwidth}
		\centering
		\includegraphics[width=\textwidth]{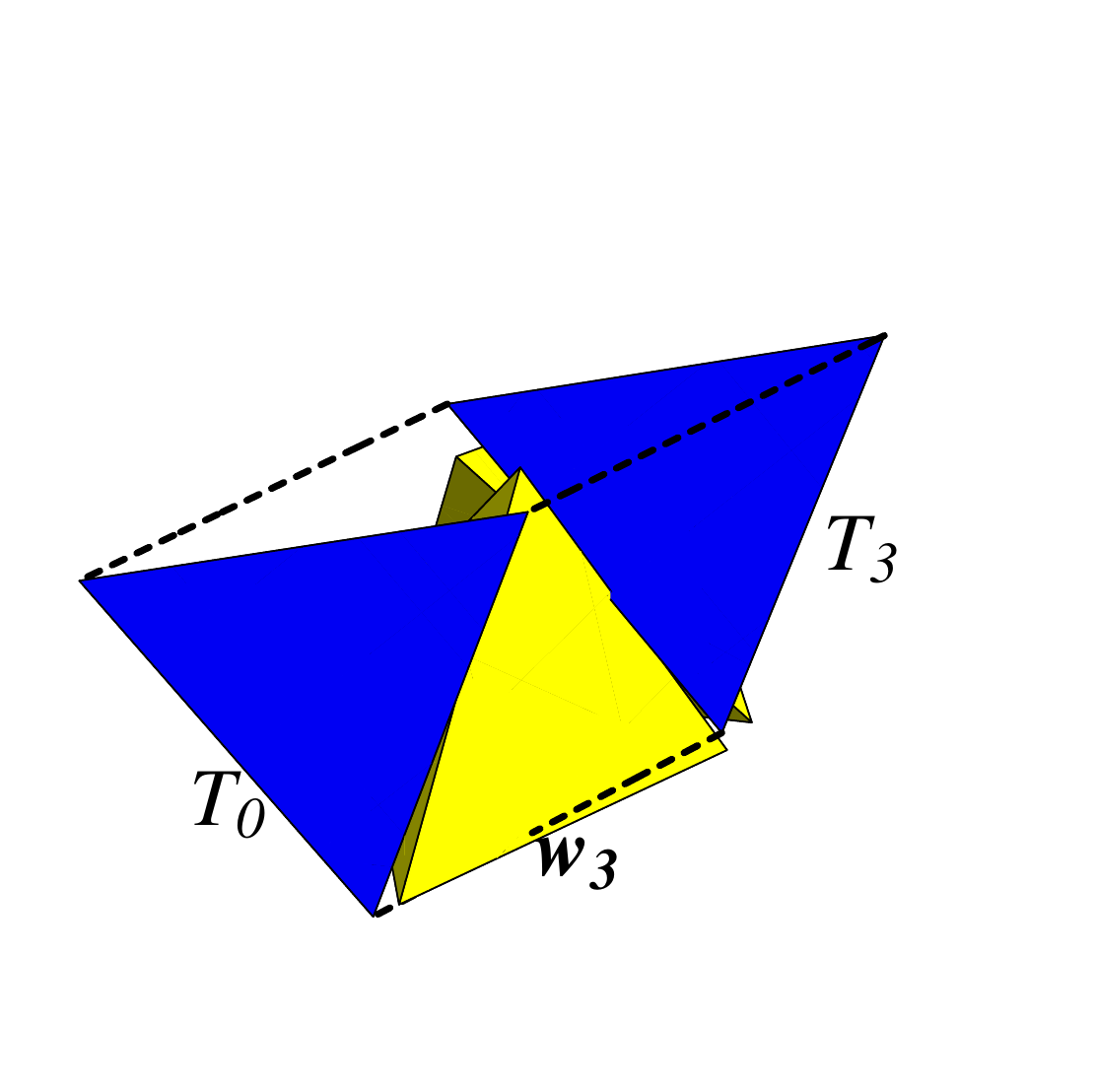}
		\caption{The vertices of $T_3$ are the vertices of $T_0$ translated by $w_3$.}
		\label{F:3philixA}
	\end{subfigure}
	\quad
	\begin{subfigure}[t]{0.45\textwidth}
		\centering
		\includegraphics[width=\textwidth]{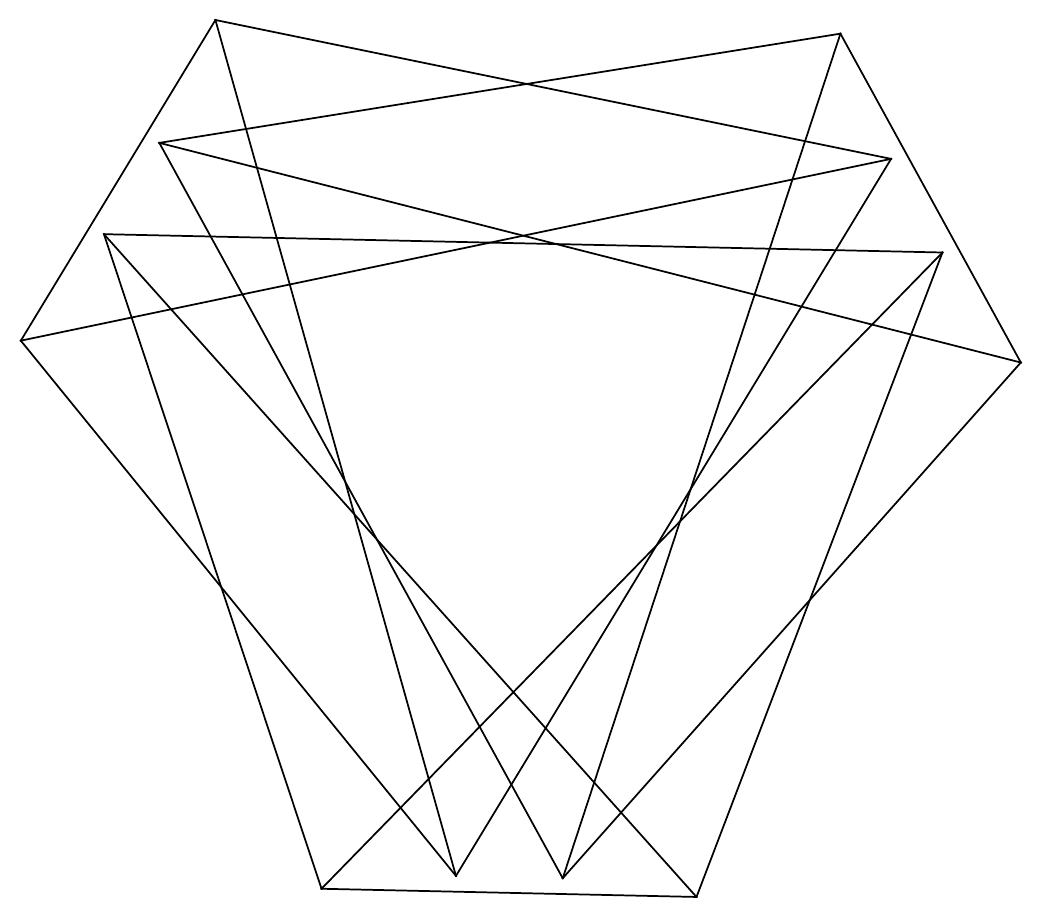}
		\caption{A projection of the 3-BC helix along its central axis.}
		\label{F:3philixB}
	\end{subfigure}
	\caption{The periodicity of the 3-BC helix.}
	\label{F:3philix}
\end{figure}

\subsection{The 3-BC helix}\label{SS:3period}

The 3-period philix (Figure~\ref{F:philicesC}) is produced here using an approach similar to that of the 5-period philix in Section~\ref{S:fperiod}. Here, a primitive set $\{ T_0, T_1, T_2 \}$ is taken such that $T_0$ is as before (see \eqref{T0} on page~\pageref{T0}) and
\newline
\begin{align}
	T_1\colon \qquad 	v_{1 0} &= \left( 0, 0, - \frac{5}{2 \sqrt{6}} \right) \label{T31}\\
					v_{1 1} &= \left( - \frac{1+3 \sqrt{5}-3 \sqrt{6-2 \sqrt{5}}}{16 \sqrt{3}}, - \frac{1+3 \sqrt{5}+ \sqrt{6-2 \sqrt{5}}}{16}, - \frac{1}{2 \sqrt{6}} \right)\notag\\
					v_{1 2} &= \left( - \frac{1+3 \sqrt{5}+3 \sqrt{6-2 \sqrt{5}}}{16 \sqrt{3}}, \frac{1+3 \sqrt{5}- \sqrt{6-2 \sqrt{5}}}{16}, - \frac{1}{2 \sqrt{6}} \right)\notag\\
					v_{1 3} &= \left( \frac{1+3 \sqrt{5}}{8 \sqrt{3}}, \frac{1}{4} \sqrt{ \frac{1}{2} \left(3-\sqrt{5} \right)}, - \frac{1}{2 \sqrt{6}} \right)\notag\\
	T_2\colon \qquad	v_{2 0} &= \left( - \frac{1}{12 \sqrt{3}}, \frac{4-\sqrt{5}}{12}, -\frac{8+3\sqrt{5}}{6\sqrt{6}} \right)\label{T32}\\
					v_{2 1} &= \left( \frac{5-3\sqrt{5}}{12\sqrt{3}}, - \frac{5+\sqrt{5}}{12}, - \frac{5}{6\sqrt{6}} \right)\notag\\
					v_{2 2} &= \left( - \frac{11+3\sqrt{5}}{24\sqrt{3}}, \frac{5+\sqrt{5}}{24}, \frac{-8+3\sqrt{5}}{6\sqrt{6}} \right)\notag\\
					v_{2 3} &= \left( - \frac{5 \left( \sqrt{3} + 3\sqrt{15} \right)}{72}, \frac{5 \left(1-\sqrt{5} \right)}{24}, - \frac{11}{6\sqrt{6}} \right).\notag
\end{align}
\newline
To generate the tetrahedra of the 3-period philix, one translates these primitive tetrahedra along an axial direction (as in Section~\ref{S:fperiod}), which now has the value
\newline
\begin{equation}\label{axisz3}
	w_3 = \left( - \frac{5+3\sqrt{5}}{12\sqrt{3}}, \frac{5-\sqrt{5}}{12}, -\frac{5}{3\sqrt{6}} \right).
\end{equation}
\newline
In other words, the 3-period philix is produced using the tetrahedra $\{ T_0, T_1, T_2 \}$ above such that
\newline
\begin{equation}\label{3periodicrelation}
	v_{\left( j + 3 k \right) i} = v_{j i} + k w_3, \text{\qquad for } k \in \mathbb{Z}.
\end{equation}
\newline

As with its 5-period sibling, the tetrahedral centroids of the 3-period philix form a helix. The corresponding values for the pitch ($p_3$) and radius ($r_3$) are substantially simpler than in the 5-period case, and are given by
\newline
\begin{align}
	p_3	&= \sqrt{\frac{5}{6}}\\
	r_3	&= \frac{\sqrt{2}}{9}
\end{align}
\newline
The corresponding parameterization to \eqref{5helixeq} is
\newline
\begin{equation}
	c \left( t \right) = r_3 \left(u_1 \cos{t} + u_2 \sin{t} \right) + \frac{t}{2\pi} w_3 + q_3,
\end{equation}
\newline
$u_1$ is as before, $u_2 = \left(\frac{1}{12} \left(\sqrt{2}-3
   \sqrt{10}\right),-\frac{1+\sqrt{5}}{2 \sqrt{6}},\frac{1}{3}\right)$, and $q_3 = \left(\frac{1}{9
   \sqrt{3}},-\frac{1}{9},-\frac{1}{9} \sqrt{\frac{2}{3}}\right)$. In this case, tetrahedral centroids lie on the helix at $t = k \frac{2 \pi}{3}, k \in \mathbb{Z}$.
   
To align the philix of this section with the \emph{z}-axis, one may use the analog of \eqref{H5transform}, with
\newline
\begin{equation}
	C_3 = \left(
\begin{array}{ccc}
 \frac{1}{12} \left(3-4 \sqrt{5}\right) & \frac{1}{12} \left(2
   \sqrt{3}+\sqrt{15}\right) & \frac{1}{12} \left(3
   \sqrt{2}+\sqrt{10}\right) \\
 \frac{1}{12} \left(2 \sqrt{3}+\sqrt{15}\right) & \frac{3}{4} &
   -\frac{-\sqrt{2}+\sqrt{10}}{4 \sqrt{3}} \\
 \frac{1}{12} \left(-3 \sqrt{2}-\sqrt{10}\right) & \frac{1}{12}
   \left(-\sqrt{6}+\sqrt{30}\right) & -\frac{\sqrt{5}}{3}
\end{array}
\right).
\end{equation}

\section{Conclusion}\label{S:conclusion}
It is known that the BC helix exhibits an aperiodic nature such that it possesses no non-trivial translational or rotational symmetries. Here we have developed modified varieties of this structure, producing helices of tetrahedra possessing both translational and rotational (in their projections) symmetries along/about their central axes. Such a structure has been designated in this writing as an \emph{m}-BC helix, and we have focused on two particular variations: the 3-BC helix (3-period philix) and the 5-BC helix (5-period philix). The construction process of these philices resembles that of the BC helix, however a rotation is added after each new tetrahedron is appended to the chain. When the value of $\beta$ given by \eqref{E:beta} is used, the relative chiralities of this rotation and the underlying chain of tetrahedra determines whether a 3-- or 5-period philix is produced.

As evidenced by Table~\ref{T:periods}, the \emph{m}-BC helix is periodic for cases other than $m = 3$ and $m = 5$. Currently, our group is working towards a general formula relating $m$ and $\beta$.
\newline
\newline
\textbf{Supplementary Information} is available in the ancillary Mathematica notebook,\newline ``\texttt{m-BC-helix-ancillary.nb}'', available with this submission at arXiv.org.

\appendix
\section{Appendix: the transformations $A_{T_k}^{f_k}$ and $B_{T_k}^{f_k}$}
The transformations $A_{T_k}^{f_k}$ and $B_{T_k}^{f_k}$ of Section~\ref{S:genphilix} have the form
\newline
\begin{align}
	A_T^f \left(v\right) &= M_T^f \left( v - c_T^f \right) + c_T^f \label{Atransform} \\
	B_T^f \left(v\right) &= R_T^f \left( v - c_T^f \right) + c_T^f \label{Btransform},
\end{align}
\newline
where $M_T^f \in \text{O}(3)$ is a reflection matrix through a mirror parallel to face $f$ of tetrahedron $T$, $R_T^f \in \text{SO}(3)$ is a rotation matrix by $\beta$ through an axis normal to the face $f$, and $c_T^f$ is the center of the tetrahedral face $f$ on $T$. The values of $M_{T_k}^{f_k}$, $R_{T_k^\prime}^{f_k}$, $T_k^\prime$, and $c_{T_k}^{f_k}$ necessary to generate the primitive tetrahedra in Sections~\ref{S:fperiod} and \ref{SS:3period} are given here.

\subsection{Transformations related to the 5-BC helix}
The reflection matrices $M_{T_0}^{f_0},\ldots,M_{T_3}^{f_3}$ are as follows:
\newline
\allowdisplaybreaks{
\begin{align}
	M_{T_0}^{f_0} &= \left(
					\begin{array}{ccc}
					 1 & 0 & 0 \\
					 0 & 1 & 0 \\
					 0 & 0 & -1
					\end{array}
					\right) \\
	M_{T_1}^{f_1} &= \left(
					\begin{array}{ccc}
					 \frac{1}{18} \left(-5-3 \sqrt{5}\right) & \frac{1}{6} \sqrt{18-\frac{14 \sqrt{5}}{3}} & -\frac{1}{9} \sqrt{23+3 \sqrt{5}} \\
					 \frac{1}{6} \sqrt{18-\frac{14 \sqrt{5}}{3}} & \frac{1}{6} \left(3+\sqrt{5}\right) & \frac{1}{3} \sqrt{1-\frac{\sqrt{5}}{3}} \\
					 -\frac{1}{9} \sqrt{23+3 \sqrt{5}} & \frac{1}{3} \sqrt{1-\frac{\sqrt{5}}{3}} & \frac{7}{9}
					\end{array}
					\right)\\
	M_{T_2}^{f_2} &= \left(
					\begin{array}{ccc}
					 \frac{1}{18} \left(11+3 \sqrt{5}\right) & -\frac{-1+\sqrt{5}}{6 \sqrt{3}} & -\frac{1}{9} \sqrt{\frac{5}{2}} \left(-3+\sqrt{5}\right) \\
					 -\frac{-1+\sqrt{5}}{6 \sqrt{3}} & \frac{1}{6} \left(3-\sqrt{5}\right) & \frac{5+\sqrt{5}}{3 \sqrt{6}} \\
					 -\frac{1}{9} \sqrt{\frac{5}{2}} \left(-3+\sqrt{5}\right) & \frac{5+\sqrt{5}}{3 \sqrt{6}} & -\frac{1}{9}
					\end{array}
					\right)\\
	M_{T_3}^{f_3} &= \left(
					\begin{array}{ccc}
					 \frac{1}{18} \left(11-3 \sqrt{5}\right) & -\frac{1+\sqrt{5}}{6 \sqrt{3}} & -\frac{1}{9} \sqrt{\frac{5}{2}} \left(3+\sqrt{5}\right) \\
					 -\frac{1+\sqrt{5}}{6 \sqrt{3}} & \frac{1}{6} \left(3+\sqrt{5}\right) & \frac{-5+\sqrt{5}}{3 \sqrt{6}} \\
					 -\frac{1}{9} \sqrt{\frac{5}{2}} \left(3+\sqrt{5}\right) & \frac{-5+\sqrt{5}}{3 \sqrt{6}} & -\frac{1}{9}
					\end{array}
					\right)
\end{align}
}
\newline
The rotation matrices $R_{T_0^\prime}^{f_0},\ldots,R_{T_3^\prime}^{f_3}$ are given by:
\newline
\allowdisplaybreaks{
\begin{align}
	R_{T_0^\prime}^{f_0} &= \left(
					\begin{array}{ccc}
					 \frac{1}{8} \left(1+3 \sqrt{5}\right) & \frac{1}{4} \sqrt{\frac{3}{2} \left(3-\sqrt{5}\right)} & 0 \\
					 -\frac{1}{4} \sqrt{\frac{3}{2} \left(3-\sqrt{5}\right)} & \frac{1}{8} \left(1+3 \sqrt{5}\right) & 0 \\
					 0 & 0 & 1
					\end{array}
					\right)\\
	R_{T_1^\prime}^{f_1} &= \left(
					\begin{array}{ccc}
					 \frac{1}{72} \left(38+15 \sqrt{5}\right) & \frac{1}{24} \sqrt{287-\frac{380 \sqrt{5}}{3}} & \frac{1}{9 \sqrt{2}} \\
					 -\frac{1}{24} \sqrt{83-\frac{104 \sqrt{5}}{3}} & \frac{1}{2}+\frac{5 \sqrt{5}}{24} & \frac{1}{6} \sqrt{14-\frac{16 \sqrt{5}}{3}} \\
					 \frac{-23+9 \sqrt{5}}{36 \sqrt{2}} & -\frac{5+\sqrt{5}}{12 \sqrt{6}} & \frac{1}{9} \left(2+3 \sqrt{5}\right)
					\end{array}
					\right)\\
	R_{T_2^\prime}^{f_2} &= \left(
					\begin{array}{ccc}
					 \frac{1}{144} \left(65+33 \sqrt{5}\right) & -\frac{-19+\sqrt{5}}{48 \sqrt{3}} & \frac{29-9 \sqrt{5}}{36 \sqrt{2}} \\
					 \frac{-41+11 \sqrt{5}}{48 \sqrt{3}} & \frac{1}{48} \left(9+17 \sqrt{5}\right) & -\frac{-1+\sqrt{5}}{12 \sqrt{6}} \\
					 \frac{11-9 \sqrt{5}}{36 \sqrt{2}} & \frac{1}{6 \sqrt{369+165 \sqrt{5}}} & \frac{1}{18} \left(11+3 \sqrt{5}\right)
					\end{array}
					\right)\\
	R_{T_3^\prime}^{f_3} &= \left(
					\begin{array}{ccc}
					 \frac{5}{36}+\frac{3 \sqrt{5}}{8} & \frac{13-2 \sqrt{5}}{24 \sqrt{3}} & \frac{-8+3 \sqrt{5}}{18 \sqrt{2}} \\
					 \frac{-17+4 \sqrt{5}}{24 \sqrt{3}} & \frac{1}{2}+\frac{5 \sqrt{5}}{24} & \frac{7-2 \sqrt{5}}{6 \sqrt{6}} \\
					 \frac{1}{36} \sqrt{83-33 \sqrt{5}} & \frac{11-7 \sqrt{5}}{12 \sqrt{6}} & \frac{1}{18} \left(11+3 \sqrt{5}\right)
					\end{array}
					\right)
\end{align}
}
\newline
The face centers $c_{T_0}^{f_0},\ldots,c_{T_3}^{f_3}$ are:
\allowdisplaybreaks{
\begin{align}
	c_{T_0}^{f_0} &= \left(0,0,-\frac{1}{2 \sqrt{6}}\right) \\
	c_{T_1}^{f_1} &= \left(-\frac{1+3 \sqrt{5}}{24 \sqrt{3}},\frac{1}{12} \sqrt{\frac{1}{2} \left(3-\sqrt{5}\right)},-\frac{7}{6 \sqrt{6}}\right) \\
	c_{T_2}^{f_2} &= \left(\frac{1}{72} \left(\sqrt{3}-7 \sqrt{15}\right),\frac{1}{24} \left(3 \sqrt{5}-1\right),-\frac{8+\sqrt{5}}{6 \sqrt{6}}\right) \\
	c_{T_3}^{f_3} &= \left(\frac{1}{72} \left(\sqrt{3}-9 \sqrt{15}\right),\frac{1}{24} \left(1+3 \sqrt{5}\right),-\frac{8+3 \sqrt{5}}{6 \sqrt{6}}\right)
\end{align}
}
\newline
The intermediate tetrahedra $T_0^\prime,\ldots,T_3^\prime$ are given by:
\newline
\allowdisplaybreaks{
\begin{align}
	T_0^\prime\colon \qquad	v_{0 0}^\prime &= \left(0,0,-\sqrt{\frac{2}{3}}-\frac{1}{2 \sqrt{6}}\right)\\
						v_{0 1}^\prime &= \left(-\frac{1}{2 \sqrt{3}},-\frac{1}{2},-\frac{1}{2 \sqrt{6}}\right)\notag\\
						v_{0 2}^\prime &= \left(-\frac{1}{2 \sqrt{3}},\frac{1}{2},-\frac{1}{2 \sqrt{6}}\right)\notag\\
						v_{0 3}^\prime &= \left(\frac{1}{\sqrt{3}},0,-\frac{1}{2 \sqrt{6}}\right)\notag\\
	T_1^\prime\colon \qquad	v_{1 0}^\prime &= \left(0,0,-\frac{5}{2 \sqrt{6}}\right)\\
						v_{1 1}^\prime &= \left(\frac{1-3 \sqrt{5}}{8 \sqrt{3}},\frac{1}{8} \left(-1-\sqrt{5}\right),-\frac{1}{2 \sqrt{6}}\right)\notag\\
						v_{1 2}^\prime &= \left(-\frac{1}{4 \sqrt{3}},\frac{\sqrt{5}}{4},-\frac{1}{2 \sqrt{6}}\right)\notag\\
						v_{1 3}^\prime &= \left(-\frac{5}{72} \left(\sqrt{3}+3 \sqrt{15}\right),\frac{5}{24} \left(\sqrt{5}-1\right),-\frac{11}{6 \sqrt{6}}\right)\notag\\
	T_2^\prime\colon \qquad	v_{2 0}^\prime &= \left(-\frac{1}{12 \sqrt{3}},\frac{1}{12} \left(\sqrt{5}-4\right),-\frac{8+3 \sqrt{5}}{6 \sqrt{6}}\right)\\
						v_{2 1}^\prime &= \left(\frac{13-11 \sqrt{5}}{24 \sqrt{3}},\frac{1}{24} \left(3+7 \sqrt{5}\right),-\frac{8+5 \sqrt{5}}{6 \sqrt{6}}\right)\notag\\
						v_{2 2}^\prime &= \left(\frac{5-3 \sqrt{5}}{12 \sqrt{3}},\frac{1}{12} \left(5+\sqrt{5}\right),-\frac{5}{6 \sqrt{6}}\right)\notag\\
						v_{2 3}^\prime &= \left(-\frac{5}{72} \left(\sqrt{3}+3 \sqrt{15}\right),\frac{5}{24} \left(\sqrt{5}-1\right),-\frac{11}{6 \sqrt{6}}\right)\notag\\
	T_3^\prime\colon \qquad	v_{3 0}^\prime &= \left(\frac{5-4 \sqrt{5}}{12 \sqrt{3}},-\frac{\sqrt{5}}{12},-\frac{11+2 \sqrt{5}}{6 \sqrt{6}}\right)\\
						v_{3 1}^\prime &= \left(\frac{13-11 \sqrt{5}}{24 \sqrt{3}},\frac{1}{24} \left(3+7 \sqrt{5}\right),-\frac{8+5 \sqrt{5}}{6 \sqrt{6}}\right)\notag\\
						v_{3 2}^\prime &= \left(-\frac{11+13 \sqrt{5}}{24 \sqrt{3}},\frac{1}{24} \left(5-\sqrt{5}\right),-\frac{8+7 \sqrt{5}}{6 \sqrt{6}}\right)\notag\\
						v_{3 3}^\prime &= \left(-\frac{5+2 \sqrt{5}}{6 \sqrt{3}},\frac{\sqrt{5}}{6},-\frac{5+2 \sqrt{5}}{6 \sqrt{6}}\right)\notag
\end{align}
}

\subsection{Transformations related to the 3-BC helix}
The reflection matrices $M_{T_0}^{f_0}$ and $M_{T_1}^{f_1}$ are as follows:
\newline
\allowdisplaybreaks{
\begin{align}
	M_{T_0}^{f_0} &= \left(
				\begin{array}{ccc}
				 1 & 0 & 0 \\
				 0 & 1 & 0 \\
				 0 & 0 & -1
				\end{array}
				\right)\\
	M_{T_1}^{f_1} &= \left(
				\begin{array}{ccc}
				 \frac{1}{18} \left(-5-3 \sqrt{5}\right) & \frac{-7+\sqrt{5}}{6 \sqrt{3}}
   & -\frac{1}{9} \sqrt{23+3 \sqrt{5}} \\
				 \frac{-7+\sqrt{5}}{6 \sqrt{3}} & \frac{1}{6} \left(3+\sqrt{5}\right) &
   -\frac{1}{3} \sqrt{1-\frac{\sqrt{5}}{3}} \\
				 -\frac{1}{9} \sqrt{23+3 \sqrt{5}} & -\frac{1}{3}
   \sqrt{1-\frac{\sqrt{5}}{3}} & \frac{7}{9}
				\end{array}
				\right)
\end{align}
}
\newline
The rotation matrices $R_{T_0^\prime}^{f_0}$ and $R_{T_1^\prime}^{f_1}$ are given by:
\allowdisplaybreaks{
\begin{align}
	R_{T_0^\prime}^{f_0} &= \left(
				\begin{array}{ccc}
				 \frac{1}{8} \left(1+3 \sqrt{5}\right) & -\frac{1}{4} \sqrt{\frac{3}{2}
   \left(3-\sqrt{5}\right)} & 0 \\
				 \frac{1}{4} \sqrt{\frac{3}{2} \left(3-\sqrt{5}\right)} & \frac{1}{8}
   \left(1+3 \sqrt{5}\right) & 0 \\
				 0 & 0 & 1
				\end{array}
				\right)\\
	R_{T_1^\prime}^{f_1} &= \left(
				\begin{array}{ccc}
				 \frac{1}{72} \left(38+15 \sqrt{5}\right) & -\frac{1}{24}
   \sqrt{287-\frac{380 \sqrt{5}}{3}} & \frac{1}{9 \sqrt{2}} \\
				 \frac{1}{24} \sqrt{83-\frac{104 \sqrt{5}}{3}} & \frac{1}{2}+\frac{5
   \sqrt{5}}{24} & -\frac{1}{6} \sqrt{14-\frac{16 \sqrt{5}}{3}} \\
				 \frac{-23+9 \sqrt{5}}{36 \sqrt{2}} & \frac{1}{12} \sqrt{\frac{5}{3}
   \left(3+\sqrt{5}\right)} & \frac{1}{9} \left(2+3 \sqrt{5}\right)
				\end{array}
				\right)
\end{align}
}
\newline
The face centers $c_{T_0}^{f_0}$ and $c_{T_1}^{f_1}$ are:
\newline
\allowdisplaybreaks{
\begin{align}
	c_{T_0}^{f_0} &= \left(0,0,-\frac{1}{2 \sqrt{6}}\right) \\
	c_{T_1}^{f_1} &= \left(-\frac{1+3 \sqrt{5}}{24 \sqrt{3}},-\frac{1}{12} \sqrt{\frac{1}{2}
   \left(3-\sqrt{5}\right)},-\frac{7}{6 \sqrt{6}}\right)
\end{align}
}
\newline
The intermediate tetrahedra $T_0^\prime$ and $T_1^\prime$ are given by:
\newline
\allowdisplaybreaks{
\begin{align}
	T_0^\prime\colon \qquad	v_{0 0}^\prime &= \left(0,0,-\sqrt{\frac{2}{3}}-\frac{1}{2 \sqrt{6}}\right)\\
						v_{0 1}^\prime &= \left(-\frac{1}{2 \sqrt{3}},-\frac{1}{2},-\frac{1}{2 \sqrt{6}}\right)\notag\\
						v_{0 2}^\prime &= \left(-\frac{1}{2 \sqrt{3}},\frac{1}{2},-\frac{1}{2 \sqrt{6}}\right)\notag\\
						v_{0 3}^\prime &= \left(\frac{1}{\sqrt{3}},0,-\frac{1}{2 \sqrt{6}}\right)\notag\\
	T_1^\prime\colon \qquad	v_{1 0}^\prime &= \left(0,0,-\frac{5}{2 \sqrt{6}}\right)\\
						v_{1 1}^\prime &= \left(-\frac{1}{4 \sqrt{3}},-\frac{\sqrt{5}}{4},-\frac{1}{2
   \sqrt{6}}\right)\notag\\
   						v_{1 2}^\prime &= \left(\frac{1-3 \sqrt{5}}{8 \sqrt{3}},\frac{1}{8}
   \left(1+\sqrt{5}\right),-\frac{1}{2 \sqrt{6}}\right)\notag\\
   						v_{1 3}^\prime &= \left(-\frac{5}{72} \left(\sqrt{3}+3 \sqrt{15}\right),-\frac{5}{24}
   \left(\sqrt{5}-1\right),-\frac{11}{6 \sqrt{6}}\right)\notag
\end{align}
}
\end{document}